\documentclass[a4paper]{article}
\usepackage[active]{srcltx}
\usepackage[cp1251]{inputenc}
\usepackage[english,ukrainian]{babel}
\newtheorem{theorem}{Теорема}
\newtheorem{lema}{Лема}

\newtheorem{guess}{Твердження}

\newtheorem{definition}{Означення}
\newtheorem{example}{Приклад}
\newtheorem{remark}{Зауваження}
\RequirePackage{latexsym} \RequirePackage{amssymb}
\RequirePackage{amscd} \RequirePackage{epsfig}
\RequirePackage{graphics} \RequirePackage{ifthen}
\newcommand{\D}{\displaystyle}

\def\<{\langle}
\def\>{\rangle}

\begin{document}
\fontsize{14pt}{16pt} \selectfont \makeatletter
\renewcommand{\@oddhead}{\hfill\Large\thepage}
\renewcommand{\@oddfoot}{}
\makeatother

УДК 517.98

М. А. Митрофанов О. В. Равський

АПРОКСИМАЦІЯ НЕПЕРЕРВНИХ ФУНКЦІЙ НА ПРОСТОРАХ ФРЕШЕ \\

У статті розглянуто апроксимацію аналітичними та $*$-аналітичними
функціями неперервних функцій на зліченно нормованих просторах
Фреше. Також знайдено критерій існування продовження неперервної
функції з всюди щільного підпростору топологічного простору на
весь простір.\\

У питанні апроксимації неперервних функцій на просторі перші
грунтовні результати були отримані ще Вейєрштрассом у  1885 році.
Проте, ці дослідження стосувалися підмножин скінченовимірних
просторів. У випадку дійсного банахового простору наступний
позитивний результат отримано Я. Курцвейлом у \cite{Ku}:
\begin {theorem} \label{1t1}.
Нехай $X$ --- сепарабельний дійсний банахів простір, що допускає
розділячий поліном, $G$ --- довільна відкрита підмножина в $X$.
Нехай $F$ --- неперервний оператор визначений в $G$ зі значеннями
в довільному банаховому просторі $Y$. Тоді для довільного
$\varepsilon > 0$ існує такий оператор $H$, аналітичний в $G$, що
нерівність:
\begin{equation}\label{eqq6_111}
||F(x)-H(x)||<\varepsilon
\end{equation}
виконується для всіх $x \in G$.
\end {theorem}
Пізніше М.~Боісо і П.~Гаєк у праці \cite{MCPH} отримали, для
апроксимації рівномірно неперервних функцій на дійсних банахових
просторах, сильніший результат при підсилені додаткових умов. У
2009 році у праці \cite{Mitr} отримано аналоги результатів Я.
Курцвейла та М.~Боісо і П.~Гаєка для комплексних банахових
просторів.

В даній роботі автори, спираючись на попередні результати отримані
для банахових просторів, досліджують можливість апроксимації
неперервних функції для деяких просторів Фреше, причому як для
дійсного, так і для комплексного випадку.

Нагадаємо основні означення. Нехай $X$ та $Y$ лінійні простори над
$\mathbb{R}$ або $\mathbb{C}$.
\begin {definition}\label{1d1}.
Відображення $B_n,$ $B_n:X^{n} \to Y$ називається
$*$-$n$-лінійним, якщо воно подається у вигляді:
$$\D B_n(x_1,...,x_k,x_{k+1},...,x_{k+m}) = \sum_{k+m=n}
c_{km}B_{km}(x_1,...,x_k,x_{k+1},...,x_{k+m}),$$ де відображення
$B_{km}:X^{k+m} \to Y$ є ненульовим $k$--лінійним відносно $x_i
\in X$, для $1 \le i \le k \in \mathbb{N}$ і $m$--антилінійним
відносно $x_{k+j} \in X$, для $1 \le j \le m \in \mathbb{N}.$
Коефіцієнти $c_{km}$ приймають значення $0$ або $1$, але принаймні
одне значення $c_{km}$ відмінне від нуля, та $n=k+m$.
\end {definition}
\begin {definition}\label{1d4}.
Відображення $F_n,$ $F_n:X \to Y$ називається $n$-однорідним
$*$-поліномом, якщо існує $*$-$n$-лінійне відображення $B_n:X^n
\to Y$ таке, що $F_n(x)=B_n(x,...,x)$, для всіх $x \in X$. У
випадку коли $n=0$, $F_0$ є тотожною константою в $Y$.
\end {definition}

\begin {definition}\label{1d5}.
Відображення $F$ з $X$ в $Y$,  називається {\em $*$-поліномом
степеня} $j$, якщо $\D F=\sum_{n=0}^j F_n,$ де $F_n$ є
$n$-однорідним $*$-поліномом  та $F_j \ne 0.$
\end {definition}

Нехай $X$ топологічний векторний простір, $Y$ нормований простір.
\begin {definition}\label{1d9}.
Відображення $F:X \to Y$ називається $*$-аналітичним, якщо для
кожної точки $x \in X$ існує окіл $V \subset X,$ $x \in V,$ такий
що $\D F(x)=\sum_{n=0}^{\infty} F_n(x),$ де $F_n$ є
$n$-однорідними неперервними $*$-поліномами і ряд $\D
\sum_{n=0}^{\infty} F_n(x)$ збігається рівномірно в околі $V$ за
нормою простору $Y.$
\end {definition}

Позначимо через $ \tilde {\mathcal H} (X,Y)$ лінійний простір всіх
$*$-аналітичних відображень з комплексного нормованого простору
$X$ в комплексний банаховий простір $Y.$

Легко бачити, що якщо $F_n$ породжуються лише формами вигляду
$B_{k0}$, то означення $F_n$ буде означенням полінома на лінійному
просторі та відповідне відображення $F$, породжене відображеннями
$F_n$, буде аналітичною функцією.

Позначимо через ${\mathcal H}(X,Y)$ лінійний простір всіх
аналітичних відображень з нормованого простору $X$ в банаховий
простір $Y.$

\begin {definition}\label{1d8}.
$*$-поліном $P:X \to \mathbb{C}$ на нормованому просторі $X$
називається розділяючим $*$-поліномом, якщо:

1.$\ P(0)=0.$

2.$\ \mid P(x) \mid \ge 1$ для кожного $x \in X$ такого, що $\Vert
x \Vert =1.$
\end {definition}

Зокрема, якщо $*$-поліном $P$ є поліномом і діє з дійсного
простору $X$ у простір $\mathbb{R}$, то він називається
розділяючим поліномом.

\begin {definition}\label{5d1}.
Нехай $X$ є комплексним нормованим простором. Будемо говорити, що
функція $Q:X \to \mathbb{C}$, $\D Q(x)=\sum_{n=1}^{\infty}Q_n(x)$
для всіх $x \in X$, де $Q_n(x)$ --- $n$-однорідні $*$-поліноми, є
рівномірно $*$-аналітичною та розділяючою, якщо вона
задовольняє наступним умовам:\\
1. Існує таке число $R_Q$, що ряд $\D \sum_{n=1}^{\infty}Q_n(x)$
збігається рівномірно в кулі радіуса
$R_Q$ з центром у довільній точці $x_0 \in X$. \\
2. Існує $\alpha \in \mathbb{R}$  таке, що множина таких $x \in
X$, що $|Q(x)|<\alpha $ є непорожньою та лежить у відкритій
одиничній кулі $B$.
\end {definition}

\begin {definition}\label{5d2}.
Нехай $X$ є дійсним нормованим простором. Будемо говорити, що
дійсна функція $Q:X \to \mathbb{R}$ є рівномірно аналітичною і
розділяючою, якщо вона задовольняє
наступним умовам:\\
1. Функція $Q$ є дійсною аналітичною на $X$ з радіусом збіжності
$R_{Q_x}$ в
кожній точці $x \in X$ більшим або рівним за $R_Q$ для деякого $R_Q>0$.\\
2. Існує $\alpha \in \mathbb{R}$ таке,  що множина таких $x \in X$
що $Q(x)<\alpha $ є непорожньою та лежить у відкритій одиничній
кулі $B$.
\end {definition}

Нам знадобляться наступні дві технічні леми. Хоча, напевне, їх
доведення є відомими, авторам не вдалося їх знайти в літературі, і
тому, заради повноти, ми наводимо доведення цих лем нижче.

Для фільтра $\mathcal F$ на топологічному просторі через
$\lim\mathcal F$ ми позначатимемо множину всіх границь фільтру
$\mathcal F$. Згідно з \cite[1.6]{ENG} ми будемо казати, що фільтр
$\mathcal F$ \em збігається до точки \em $x$, якщо $x\in\lim
\mathcal F$, та що фільтр $\mathcal F$ є \em збіжним, \em якщо
$\lim\mathcal F\not=\emptyset$.

Нехай $X$ --- топологічний простір, $Y$ --- регулярний
топологічний простір, $D$ -- щільна підмножина простору $X$, $g:D
\to Y$ відображення та $\mathcal F$ --- фільтр на просторі $X$.
Через $g(\mathcal F)$ позначимо фільтр на просторі $Y$, породжений
базою  $\{g(F):F\in \mathcal F\}$. Для кожної точки $x\in X$ через
$\mathcal F_x$ позначимо слід фільтру всіх околів точки $x$ на
множині $D$.

\begin {lema}\label{TopConv}. Неперервне відображення $f:D \to Y$ продовжується
до неперервного відображення $\widetilde{f}:X \to Y$ тоді і лише
тоді, коли для довільної точки $x\in X$ фільтр $f(\mathcal F_x)$ є
збіжним.
\end {lema}
{\sc Доведення.} \em Необхідність. \em Нехай $x\in X$ --- довільна
точка. Тоді, оскільки відображення $f$ є неперервним, за
твердженням \cite[1.6.10]{ENG}, виконується включення
$$\widetilde{f}(x)\in\widetilde{f}(\lim\mathcal F_x)\subset
\lim\widetilde{f}(\mathcal F_x)=\lim f(\mathcal F_x),$$ отже
фільтр $f(\mathcal F_x)$ є збіжним.

\em Достатність. \em Оскільки простір $Y$ є хаусдорфовим, то, за
твердженням \cite[1.6.11]{ENG}, для кожної точки $x\in X$ фільтр
$f(\mathcal F_x)$ має єдину границю. Побудуємо відображення
$\widetilde{f}$, прийнявши $\{ \widetilde{f}(x) \} =\lim
f(\mathcal F_x)$, для кожної точки $x\in X$. За неперервністю $f$
на $D$ маємо $\{\widetilde{f}(x)\}=\lim f(\mathcal F_x)\supset
f(\lim \mathcal F_x)=\{f(x)\}$. Отже, відображення $\widetilde{f}$
є продовженням відображення $f$.

Покажемо тепер неперервність відображення $\widetilde{f}$. Нехай
$\widetilde{V}$ --- довільна непорожня відкрита підмножина
простору $Y$ та $x\in \widetilde{f}^{-1}(\widetilde{V})$. За
регулярністю простору $Y$, існує відкрита підмножина
$\widetilde{U}$ така, що $\widetilde{f}(x)\in
\widetilde{U}\subset\overline{\widetilde{U}}\subset
\widetilde{V}$. Оскільки $\lim f(\mathcal F_x)\subset
\widetilde{U}$, то існує відкритий окіл $U$ точки $x$, такий, що
$f(U\cap D)\subset\widetilde{U}.$ Тоді для кожної точки $x'\in U$
маємо $\{ \widetilde{f}(x') \} =\lim f(\mathcal F_{x'})\subset
\overline{\widetilde{U}}\subset\widetilde{V}$.

\ $\square$

Топологічний простір $X$ називається \em {простором Фреше-Урисона}
\em (див. наприклад \cite[1.6]{ENG}), якщо для довільної $A\subset
X$ та довільної $x\in\overline A$ існує послідовність $\{x_n\}$
точок множини $A$, збіжна до $x$. Кожен простір з першою аксіомою
зліченності (а, отже, і кожен метризовний простір) є простором
Фреше-Урисона \cite[1.6.14]{ENG}.

\begin {lema}\label{8}. Нехай $X$ --- простір Фреше-Урисона, $Y$ ---
регулярний топологічний простір, $D$ --- щільна підмножина
простору $X$. Неперервне відображення $f:D \to Y$ продовжується до
неперервного відображення $\widetilde{f}:X \to Y$ тоді і лише
тоді, коли для довільної збіжної в $X$ послідовності $\{x_n\}$
точок з $D$ послідовність $\{f(x_n)\}$ теж є збіжною.
\end {lema}
{\sc Доведення.} \em Необхідність. \em Нехай $x\in X$ --- довільна
точка, та $\{x_n\}\subset D$ --- збіжна до $x$ послідовність.
Позначимо через $\mathcal S_x$ фільтр на $D$, породжений базою
$\{\{x_n:n\ge m\}:m\in \mathbb{N} \}$. За Лемою \ref{TopConv}
фільтр $f(\mathcal F_x)$ є збіжним. Легко бачити, що фільтр
$\mathcal S_x$ є тонкішим за фільтр $\mathcal F_x$. Оскільки
фільтр $f(\mathcal F_x)$ є збіжним, то за твердженням
\cite[1.6.8]{ENG}, фільтр $f(\mathcal S_x)$ також є збіжним, і,
отже, послідовність $\{f(x_n)\}$ теж є збіжною.

\em Достатність. \em Нехай $x\in X$ --- довільна точка, та
$\{x_n\}\subset D$ --- збіжна до $x$ послідовність. Нехай
послідовність $\{f(x_n)\}$ збігається до точки $y\in Y$. Ми
покажемо, що фільтр $f(\mathcal F_x)$ теж збігається до точки $y$.
Дійсно, припустимо протилежне. Тоді існує такий окіл
$\widetilde{V}$ точки $y$, що для довільного околу $U$ точки $x$
існує точка $x_U\in U\cap D$ така, що $f(x_U)\in Y\backslash
\widetilde{V}$. Позначимо сім'ю всіх околів точки $x$ через
$\mathcal N_x$ і приймемо $A=\{x_U:U \in \mathcal N_x \}$. З
побудови $A$ випливає, що $x\in\overline{A}$, тому існує
послідовність точок $\{x_n'\} \subset A$ збіжна до $x$. Задамо
послідовність  $\{x_n''\}\subset D$ наступним чином, поклавши
$x_{2n}''=x_n$ та $x_{2n-1}''=x_n'$ для всіх $n \in \mathbb{N}$.
Тоді послідовність $\{x_n''\}$ є збіжною до $x$, тому, за умовою
леми, послідовність $\{f(x_n'')\}$ теж є збіжною до деякої точки
$y''\in Y$. Зрозуміло, що послідовність $\{f(x_n)\}$
 є підпослідовністю послідовності $\{f(x_n'')\}$. Оскільки
$\{f(x_n'')\}$ збігається до $y''$, а $\{f(x_n)\}$ збігається до
$y$, то $y''=y$. Тому послідовність $\{f(x_n')\}$ теж збігається
до точки $y$. Але, за побудовою множини $A$, $\{f(x_n')\}\subset
Y\backslash \widetilde{V}$, отже $y\in Y\backslash \widetilde{V}$.
Це суперечить тому, що $y\in \widetilde{V}$. Таким чином фільтр
$f(\mathcal F_x)$ збігається до точки $y$. Тому, за лемою
\ref{TopConv}, відображення $f:D \to Y$ продовжується до
неперервного відображення $\widetilde{f}:X \to Y$.

\ $\square$

Топологічний простір $X$ називається \em {секвенціальним
простором}, \em \cite[1.6]{ENG} якщо множина $A\subset X$ замкнена
тоді і тільки тоді, коли разом з кожною послідовністю вона містить
всі її границі. Кожен простір Фреше-Урисона є секвенціальним.
\cite[1.6.14]{ENG}

Наступний приклад показує, що в лемі \ref{8t30} умову ''$X$ ---
простір Фреше-Урисона'' не можна послабити до умови ''$X$ ---
секвенціальний простір''.

\begin{example} Нехай $X$ --- секвенціальний досконало нормальний
простір з прикладу \cite[1.6.19]{ENG}. Тобто,
$X=\{0\}\cup\bigcup_{i=1}^\infty X_i$, де
$X_i=\{1/i\}\cup\bigcup_{j=i^2}^\infty \left\{\frac 1i+\frac
1j\right\}$. Тоді $X_i\cap X_k=\emptyset$ при $i\not=k$. Топологія
на $X$ породжуеться наступною системою околів. Всі точки $\frac
1i+\frac1j$ --- ізольовані точки простору $X$. Для точок вигляду
$\frac 1i$ візьмемо в якості сім'ї околів сім'ю всіх множин
$X_i\backslash \bigcup_{j=i^2}^k \left\{\frac 1i+\frac 1j\right\}$
для $k=i^2,i^2+1,\dots.$ Нарешті, в якості елементів бази в точці
$0$ візьмемо всі множини, отримані з $X$ викиданням скінченної
кількості членів $X_i$ і скінченної кількості точок вигляду $\frac
1i+\frac 1j$ у всіх $X_i$, котрі залишились.

Приймемо $D=X\backslash (\{0\}\cup \bigcup_{i=1}^\infty\{1/i\})$ і
$Y=X\backslash\{0\}$ та розглянемо відображення $f:D\to Y$ таке,
що $f(x)=x$ для всіх точок $x\in D$.

Нехай $\{x_n\}$ --- довільна з послідовність точок з $D$, збіжна
до точки $x\in X$. Легко показати, що $x\not=0$. Оскільки
відображення $f$ є тотожнім на просторі $D$, то послідовність
$\{f(x_n)\}$ теж збіжна до точки $x$ в просторі $Y$.

Припустимо тепер, що існує неперервне продовження
$\widetilde{f}:X\to Y$ відображення $f$. Зафіксуємо довільне число
$i\in\mathbb{N}$. Послідовність $\{\frac 1i+\frac
1j:j\in\mathbb{N}, j\ge i^2\}$ збігається до точки $1/i\in X$.
Тоді, за лемою \cite[1.6.15]{ENG}, послідовність
$\{\widetilde{f}(\frac 1i+\frac 1j):j\in\mathbb{N}, j\ge i^2\}$
збігається до точки $\widetilde{f}(1/i)$. Але, оскільки ця
послідовність має єдину границю в просторі $Y$, то
$\widetilde{f}(1/i)=1/i$. Застосуємо лему \cite[1.6.15]{ENG} до
послідовності $\{1/i\} \subset X$. Отримаємо, що
$\lim\{\widetilde{f}(1/i)\}\supset \widetilde{f}(\lim\{1/i\})=\{
\widetilde{f}(0) \}$. Але це неможливо, бо послідовність $\{1/i\}$
не є збіжною в просторі $Y$. Ця суперечність показує, що не існує
неперервного продовження відображення $f$ на простір $X$.
\end{example}

Наступна лема доведена у праці \cite[4.3.17]{ENG}.

\begin{lema}\label{8c12} Нехай $(X,\rho)$ --- метричний простір, $(Y,\sigma)$
--- повний метричний простір, $D$ --- щільна підмножина простору $X$. Тоді кожне
відображення $f:D \to Y$, рівномірно неперервне відносно $\rho$ та
$\sigma$, продовжується до відображення $\widetilde{f}:X \to Y$,
рівномірно неперервного відносно $\rho$ та $\sigma$.
\end{lema}

Зауважимо, що існують неперервні функції, що задовольняють умовам
леми \ref{8}, але не є рівномірно неперервними. Наприклад, функція
$f(x)=\frac{1}{x}$, що діє зі всюди щільної підмножини
$\mathbb{R}_+$ в $\mathbb{R}_+$ продовжується на все
$\mathbb{R}_+$ та задовольняє умовам леми \ref{8}, але не є
рівномірно неперервною.

Нагадаємо, що лінійний топологічний простір $X$ є простором Фреше,
якщо $X$ є метризовним повною метрикою локально опуклим простором.
Відомо, що це означення еквівалентно до наявності на $X$ зліченної
системи напівнорм $\{ p_n \}_{n \in \mathbb{N}}$, які задають
повну метрику $\rho$ на $X$ таким способом \cite{SHEFER}:
\begin{equation}\label{10}
\rho(x,y)=\sum_n \frac{1}{2^n} \frac{p_n(x-y)}{1+p_n(x-y)}.
\end{equation}

Топологія, що породжується метрикою $\rho$ на просторі $X$, є
найслабшою топологією відносно якої всі напівнорми $p_n$ є
неперервними. Базу околів нуля цієї топології утворює сім'я
$\{U_{p_n}(0):n\in\mathbb{N}\}$, де $U_{p_n}(0)= \{ x \in X: \
p_k(x)< \frac{1}{n} \ \mbox{ для всіх } \ 1 \leqslant k \leqslant
n \}$. При цьому послідовність $\{x_n:n\in\mathbb{N}\}$ точок
простору $X$ прямує до точки $x_0\in X$ тоді і лише тоді, коли для
довільного $k \in \mathbb{N}$, послідовність
$\{p_k(x_n-x_0):n\in\mathbb{N}\}$ прямує до нуля, коли $n$ прямує
до нескінченності. Для нормованого простору $Y$ функція $f:X \to
Y$ є неперервною в точці $x_0$, якщо для довільної послідовності
$\{ x_n \} \subset X$, що прямує до $x_0$ послідовність $\{f(x_n)
\} \subset Y$ прямує до $f(x_0)$.

Наступне твердження є доведеним у літературі.

\begin {guess} \label{8t26}.
Нехай простір $X$ є простором Фреше з системою напівнорм $\{ p_n
\}_{n \in \mathbb{N}}$, а $Y$ банаховим простором. Функція $f:X
\to Y$ є неперервною в точці $x_0 \in X$, якщо існує $k \in
\mathbb{N}$ таке, що $f$ неперервна відносно напівнорми $p_k$.
\end {guess}

Зафіксуємо напівнорму $p_n$ на $X$. Відомо, що $\ker p_n$ є
замкненим лінійним підпростором. Позначимо через $\widetilde{X_n}$
поповнення фактор простору $X/\ker p_n$. У випадку коли $p_n$ є
нормами $X=X/\ker p_n$ та фактор норма співпадає з $p_n$.

Припустимо, що простір $X$ є сепарабельним та для довільного $X_n$
існує розділяючий поліном (розділяюча рівномірно аналітична
функція). Чи буде випливати з цих умов, що кожна неперервна
(рівномірно неперервна) функція апроксимується аналітичними на
$X$? Нижче ми даємо часткову відповідь на ці питання.

\begin {theorem} \label{8t27}.
Нехай простір $X$ є сепарабельним дійсним простором Фреше зі
зліченною системою норм $\{ p_n \}_{n \in \mathbb{N}}$, а $Y$ ---
банаховим простором. Нехай для довільного $n \in \mathbb{N}$
простір $X_n=(X,p_n)$ допускає розділяючий поліном. Тоді кожна
функція $f:X \to Y$, для якої існує $k \in \mathbb{N}$, що для
довільної фундаментальної послідовності точок $\{x_n\}$ в $X_k$,
послідовність $\{f(x_n)\} \subset Y$ є збіжною, наближається
аналітичними рівномірно на всьому $X$.
\end {theorem}
{\sc Доведення.} Оскільки для довільної фундаментальної
послідовності точок $\{x_n\}$ в $X_k$, послідовність $\{f(x_n)\}
\subset Y$ є збіжною, то функція $f$ є неперервною на просторі
$X_k$, а отже за твердженням \ref{8t26} і на просторі $X$.
Зауважимо, що якщо простір $X_k=(X,p_k)$ не є банаховим простором,
то він неповний відносно норми $p_k$. Оскільки $p_n$ є системою
норм, то для довільного $n \in \mathbb{N}$ носії просторів $X_n$
та $X$ співпадають. Поповнимо простір $X_k$ відносно норми $p_k$
до банахового простору $\widetilde{X_k}$. Тоді, оскільки $X_k$ є
щільним в $\widetilde{X_k}$, і $X_k$ допускає розділяючий поліном,
то $\widetilde{X_k}$ теж допускає розділяючий поліном. З умов
теореми, щільності $X_k$ в $X$ та за лемою \ref{8}, існує
неперервне продовження $\widetilde{f}$ відображення $f$ на простір
$\widetilde{X_k}$. За теоремою Курцвейла \ref{1t1} функція
$\widetilde{f}$ рівномірно наближається послідовністю аналітичних
функцій $\{ \widetilde{f_m} \}$ на просторі $\widetilde{X_k}$.
Звуження $f_m$ функції $\widetilde{f_m}$ на простір $X_k$ для
довільного $m \in \mathbb{N}$ є аналітичним за означенням. Тому
функція $f$ рівномірно наближається послідовністю $\{ f_m \}$
аналітичних функцій на просторі $X_k$. Легко бачити, що
відображення $f_m$ є аналітичним на просторі $X$, а отже функція
$f$ рівномірно наближається послідовністю $\{ f_m \}$ аналітичних
функцій на просторі $X$.

\ $\square$

\begin {theorem} \label{8t28}.
Нехай простір $X$ є сепарабельним дійсним простором Фреше зі
зліченною системою норм $\{ p_n \}_{n \in \mathbb{N}}$, а $Y$
--- банаховим простором. Нехай для довільного $n \in \mathbb{N},$
простір $X_n=(X,p_n)$ допускає рівномірно аналітичну і розділяючу
функцію. Тоді кожна рівномірно неперервна функція $f:X \to Y$
наближається аналітичними рівномірно на всьому $X$, якщо існує $k
\in \mathbb{N}$, що $f$ є рівномірно неперервною на $X_k$.
\end {theorem}

Доведення цієї теореми є аналогічним до доведення попередньої
теореми, тільки замість леми \ref{8}, ми використовуємо лему
\ref{8c12} та замість теореми Курцвейла \ref{1t1} ми
використовуємо основний результат М.~Боісо і П.~Гаєка
\cite[теорему 1, стор. 83]{MCPH}.

\ $\square$

Зауважимо, що у теоремі \ref{8t28} вимогу рівномірної
неперервності $f$ можна опустити.

Аналогічно, спираючись на леми \ref{8}, \ref{8c12} та результати
роботи \cite{Mitr} можна довести дві наступні теореми.

\begin {theorem} \label{8t29}.
Нехай простір $X$ є сепарабельним комплексним простором Фреше зі
зліченною системою норм $\{ p_n \}_{n \in \mathbb{N}}$, а $Y$
--- банаховим простором. Нехай для довільного $n \in \mathbb{N},$
$X_n=(X,p_n)$ допускає розділяючий $*$-поліном. Тоді кожна функція
$f:X \to Y$ така, що існує $k \in \mathbb{N}$, таке що для
довільної фундаментальної послідовності точок $\{x_n\}$ в $X_k$,
послідовність $\{f(x_n)\} \subset Y$ є збіжною наближається
$*$-аналітичними рівномірно на всьому $X$.
\end {theorem}

\begin {theorem} \label{8t30}.
Нехай простір $X$ є сепарабельним комплексним простором Фреше, зі
зліченною системою норм $\{ p_n \}_{n \in \mathbb{N}}$, а $Y$
--- банаховим простором. Нехай для довільного $n \in \mathbb{N},$
$X_n=(X,p_n)$ допускає рівномірно $*$-аналітичну та розділяючу
функцію. Тоді кожна рівномірно неперервна функція $f:X \to Y$
наближається $*$-аналітичними рівномірно на всьому $X$, якщо існує
$k \in \mathbb{N}$, таке що $f$ є рівномірно неперервною на $X_k$.
\end {theorem}

Результати теорем \ref{8t29} (\ref{8t27}) можна узагальнити
наступним чином.

\begin{remark}
Нехай простори $X$ та $Y$ такі як у теоремах \ref{8t29}
(\ref{8t27}). Через $C(X,Y)$ ми позначимо простір неперервних
функцій з простору $X$ у простір $Y$, наділений топологією
рівномірної збіжності. Замикання $\overline{\tilde {\mathcal H}
(X,Y)}$ ($\overline{\mathcal H (X,Y)}$)простору $ \tilde {\mathcal
H} (X,Y)$ ($ \mathcal H (X,Y)$) у просторі $C(X,Y)$ --- це
лінійний простір всіх неперервних функцій з простору $X$ у простір
$Y$, які є рівномірно апроксимовними $*$-аналітичними
(аналітичними) функціями з простору $X$ у простір $Y$.

Нехай для всіх $k \in \mathbb{N}$ на просторі $X_k=(X,p_k)$ існує
розділяючий поліном. Зафіксуємо довільний індекс $k \in
\mathbb{N}$. Через $\widetilde{C}(X_k,Y)\subset C(X,Y)$ ми
позначимо лінійний простір неперервних функцій з простору $X_k$ у
простір $Y$, які продовжуться до неперервних функцій з поповнення
$\widetilde{X_k}$ нормованого простору $X_k$ у простір $Y$. За
теоремами \ref{8t29} (\ref{8t27}) $\widetilde{C}(X_k,Y) \subset
\overline{\tilde {\mathcal H} (X,Y)}$ ($\widetilde{C}(X_k,Y)
\subset \overline{\mathcal H (X,Y)}$). Через $\overline{C}(X,Y)$
позначимо замикання лінійного підпростору, породженого множиною
$\bigcup\{\widetilde{C}(X_k,Y):k\in \mathbb{N}\} $ у просторі
$C(X,Y)$. Тоді $\overline{C}(X,Y) \subset \overline{\tilde
{\mathcal H} (X,Y)}$ ($\overline{C}(X,Y) \subset
\overline{\mathcal H (X,Y)}$).
\end{remark}

До наведеного у прикладі 4 праці \cite[3.5.3]{Kl-Fm}
злічено-гільбертового простору швидкоспадних послідовностей з
системою норм $\{||x||_k:k\in\mathbb{N}\}$, де $||x||_k=\left(
\sum_{n=1}^{\infty}n^k x_n^2 \right)^{\frac{1}{2}}$, ми можемо
застосувати попередні теореми. Оскільки цей простір не є
банаховим, то в роботі знайдено нові простори, на яких неперервні
функції певного вигляду допускають рівномірну апроксимацію
аналітичними.

Наступний приклад показує нетривіальність лінійного простору
$\overline{C}(X,Y)$.

\begin {example} \label{1e41}.
Нехай $X$ --- простір Фреше, що не є банаховим, з топологією,
заданою системою норм $\{ p_n \}$. Нехай $f(x)=\rho(x,0)$ для всіх
$x \in X$ де $\rho$ --- метрика на просторі $X$, визначена
формулою (\ref{10}). За побудовою, $f \in
\overline{C}(X,\mathbb{R})$. Припустимо, що існує індекс $k \in
\mathbb{N}$ такий, що $f \in \overline{C}(X_k,\mathbb{R})$.
Оскільки простір $X$ не є банаховим, топологія Фреше на $X$ строго
сильніша за топологію задану нормою $p_k$. Отже існує таке число
$\varepsilon > 0$, що множина $\{ x \in X\ : \ \rho(x,0)<
\varepsilon \}=f^{-1}(-\infty,\varepsilon)$ не містить жодного
відкритого в $X_k$ околу нуля, що суперечить неперервності функції
$f$ на просторі $X_k$.
\end {example}

Наведемо приклад, який показує суттєвість умови локальної
опуклості для апроксимації неперервних функції.

\begin{example} \label{8e6}
Нехай $X$ --- це простір $L_p[0,1]$ де $0<p<1$. Метрика на $X$
задається у наступний спосіб:
$$
\rho(x,y)=\left( \int_0^1 |f(t)-g(t)|^p dt\right)^{1/p}.
$$
Відомо, що відносно метрики $\rho$ простір $X$ є повним, але не
локально опуклим. Зокрема на цьому просторі не існує жодного
ненульового лінійного неперервного функціоналу (див.
\cite{SHEFER}, стор. 44, 49). Отже на $X$ не існує жодної
аналітичної функції відмінної від сталої (бо похідна Фреше
аналітичної функції є лінійним неперервним функціоналом). З іншого
боку, на $X$
 існують неперервні не сталі функції, наприклад $F(f)=\rho(0,f)$,
 які, отже, не наближаються аналітичними.

\end{example}

Автори висловлюють вдячність Тарасу Банаху та Андрію Загороднюку,
за участь у обговоренні статті.

М. А. Митрофанов А. В. Равский

АППРОКСИМАЦИЯ НЕПРЕРЫВНЫХ ФУНКЦИЙ НА ПРОСТРАНСТВАХ ФРЕШЕ

В статье рассмотрена аппроксимация аналитическими и
$*$-аналитическими функциями непрерывных функций на счетно
нормированых пространствах Фреше. Найден критерий существования
продолжения непрерывной функции с всюду плотного подпространства
топологического пространства на всё пространство.

M. A. Mytrofanov A. V. Ravsky

APPROXIMATION OF CONTINUOUS FUNCTIONS ON FRECHET SPACES

We consider approximations of a continuous function on a countable
normed Fr\'{e}chet space by analytic and $*$-analytic. Also we
found a criterium of the existence of an extension of a continuous
function from a dense subspace of a topological space onto the
space.

Ін-т прикл. проблем механіки і математики\\
ім. Я. С. Підстригача НАН України, Львів

\end{document}